\documentclass{amsart}
\usepackage{amssymb,latexsym}
\usepackage{amsmath}
\usepackage{amsthm}

\usepackage{graphicx}

\newtheorem{thm}{Theorem}[section]
\newtheorem{prop}[thm]{Proposition}

\newtheorem{cor}[thm]{Corollary}
\newtheorem{rem}[thm]{Remark}
\newtheorem{ex}[thm]{Example}

\begin{document}
\title{Face numbers of cubical barycentric subdivisions}
\author{Christina~Savvidou}
\address{Department of Mathematics\\
University of Athens\\
Panepistimioupolis\\
15784 Athens, Greece}
\email{savvtina@math.uoa.gr}

\thanks{Supported by the Cyprus State Scholarship Foundation.}

\begin{abstract}
The cubical barycentric subdivision $\mathrm{sd}_c(K)$ of a cubical
complex $K$ is introduced as an analogue of the barycentric
subdivision of a simplicial complex. Explicit formulas for the short
and long cubical $h$-vector of $\mathrm{sd}_c(K)$ are given, in
terms of those of $K$. It is deduced that symmetry and nonnegativity
of these $h$-vectors, as well as real rootedness of the short
cubical $h$-polynomial, are preserved under cubical barycentric
subdivision. The asymptotic behavior of the short and long cubical
$h$-vectors of successive cubical barycentric subdivisions of $K$ is
also determined.
\end{abstract}

\maketitle

\section{Introduction}
The present work is partly motivated by \cite{brentiwelker}. In that
article, Brenti and Welker study the transformation of the
$h$-vector of a simplicial complex $\Delta$ under barycentric
subdivision. They express the entries of the $h$-vector of the
barycentric subdivision of $\Delta$ as nonnegative integer linear
combinations of those of the $h$-vector of $\Delta$. In particular,
they show that symmetry and nonnegativity of the $h$-vector are
preserved under barycentric subdivision. Moreover, they prove that
if the $h$-vector of $\Delta$ is nonnegative, then the
$h$-polynomial of the barycentric subdivision of $\Delta$ has only
real roots.

It is natural to inquire whether similar results hold for
non-simplicial complexes. We study the transformation of the short
and the long cubical $h$-vector of a cubical complex $K$ under a
natural cubical analogue of simplicial barycentric subdivision,
which we call \emph{cubical barycentric subdivision}. The cubical
barycentric subdivision $\mathrm{sd}_c(K)$ is a cubical complex
which subdivides $K$, so that the poset of nonempty faces of
$\mathrm{sd}_c(K)$ is isomorphic to the set of closed intervals in
the poset of nonempty faces of $K$, partially ordered by inclusion
(see Section 2.3 for a precise definition). Our main result
expresses the entries of the short and long cubical $h$-vector of
$\mathrm{sd}_c(K)$ explicitly as nonnegative linear combinations of
the entries of the corresponding $h$-vector of $K$ (Theorems
\ref{hsc} and \ref{thm h^c}). From these expressions we deduce that
symmetry and nonnegativity of the short and long cubical $h$-vector,
as well as real rootedness of the short cubical $h$-polynomial, are
preserved under cubical barycentric subdivision. We also study the
asymptotic behavior of the short and long cubical $h$-polynomials
under successive cubical barycentric subdivisions (Corollaries
\ref{sch-polynomial 4 n-th subdivision} and \ref{ch-polynomial 4
n-th subdivision}).

\textsc{Acknowledgement.} I would like to thank Christos
Athanasiadis for suggesting this problem and for his extensive
comments on earlier versions of this paper and Volkmar Welker for
useful discussions.

\section{Preliminaries} \label{preliminaries}
This section reviews basic definitions concerning cubical complexes
and their $f$-vectors and $h$-vectors and introduces the concept of
cubical barycentric subdivision for such complexes. For backround
and any undefined terminology on partially ordered sets and on
polytopes and polyhedral complexes we refer the reader to
\cite[Chapter 3]{stanley-enum} and  to \cite{grunbaum,ziegler},
respectively.

\subsection{Cubical complexes}
We denote by $C_d$ the standard $d$-dimensional cube
$[0,1]^d\subseteq \mathbb{R}^d$. Any polytope which is
combinatorially isomorphic to $C_d$ is said to be a
\emph{combinatorial $d$-cube}. A \emph{cubical complex} is a finite
collection $K$ of combinatorial cubes in $\mathbb{R}^n$, such that
(i) every face of an element of $K$ also belongs to $K$ and (ii) the
intersection of any two elements of $K$ is a face of both. The
elements of $K$ are called \emph{faces}.

We denote by $\mathcal{F}(K)$ the face poset of $K$, meaning the set
of faces of $K$, partially ordered by inclusion. This poset is a
meet-semilattice. The empty face is the minimum element of
$\mathcal{F}(K)$ and the vertices of $K$ are its atoms. The maximal
elements of $\mathcal{F}(K)$ are called \emph{facets}. The dimension
of $K$, denoted by $\dim(K)$, is defined as the maximum dimension of
a face.

\subsection{Face enumeration} Let $K$ be a $(d-1)$-dimensional cubical complex. We denote by $f_i(K)$
the number of $i$-dimensional faces of $K$. The \emph{$f$-vector} of
$K$ is defined as
\[f(K)=(f_0(K),f_1(K),\ldots,f_{d-1}(K)).\] The polynomial
\[f_K(x) \; = \; \sum_{i=0}^{d-1}f_i(K)\, x^i\] is
called the \emph{$f$-polynomial} of $K$. The \emph{short cubical
$h$-polynomial} of $K$ is defined in \cite{adin} by the equation
\begin{equation}\label{h^sc def}
h^{(sc)}_K(x) \; = \; \sum_{i=0}^{d-1}h^{(sc)}_i(K)\; x^i \; = \;
\sum_{j=0}^{d-1} f_j(K) (2x)^j (1-x)^{d-1-j}.
\end{equation}
The vector
\[h^{(sc)}(K) \; = \; (h^{(sc)}_0(K),h^{(sc)}_1(K),
\ldots,h^{(sc)}_{d-1}(K))\] of coefficients of this polynomial is
called the \emph{short cubical $h$-vector} of $K$. The polynomials
$f_K(x)$ and $h_K^{(sc)}(x)$ are related by the equations
\begin{equation}\label{syndesh polyn. short h me f}
h_K^{(sc)}(x)\; = \; (1-x)^{d-1} f_K\left(\frac{2x}{1-x}\right)
\end{equation}
and
\begin{equation}\label{syndesh polyn. f me short h}
2^{d-1} f_K(x)\; = \; (x+2)^{d-1}
h^{(sc)}_K\left(\frac{x}{x+2}\right).
\end{equation}
As a result, the entries of $f(K)$ can be expressed in terms of
those of $h^{(sc)}(K)$ and vice versa by the equations
\begin{equation}\label{f synarthsei h^sc}
f_j(K)\; = \;2^{-j}\sum_{i=0}^j\binom{d-1-i}{d-1-j}h^{(sc)}_i(K)
\end{equation}
and
\begin{equation}\label{h^sc synarthsei f}
h^{(sc)}_i(K)\; = \;\sum_{j=0}^i
\binom{d-1-j}{d-1-i}(-1)^{i-j}2^jf_j(K).
\end{equation}
The (long) \emph{cubical $h$-vector}
\[h^{(c)}(K)=(h^{(c)}_0(K),h^{(c)}_1(K),\ldots, h^{(c)}_{d}(K))\]
of $K$ is defined in \cite{adin} by the recursive formula
\begin{equation}\label{h^sc synarthsei h^c}
h^{(sc)}_i(K) \; = \; h^{(c)}_i(K) + h^{(c)}_{i+1}(K), \; \textrm{
for } 0\leq i\leq d-1
\end{equation}
and the initial condition $h^{(c)}_0(K)=2^{d-1}$. The polynomial
\[
h^{(c)}_K(x)\; = \;\sum_{i=0}^d h^{(c)}_i(K)\, x^i
\]
is called the (long) \emph{cubical $h$-polynomial} of $K$. The short
and long cubical $h$-polynomials of $K$ are related by the
equation
\begin{equation}\label{rel short-long h-pol}
(1+x) h^{(c)}_K(x) \; = \; 2^{d-1} + x h^{(sc)}_K(x) + 2^{d-1}
(-x)^{d+1} \tilde{\chi}(K),
\end{equation}
where
\[
\tilde{\chi}(K) \;= \; -1 +\sum_{i=0}^{d-1}(-1)^if_i(K) = -1+f_K(-1)
\]
is the reduced Euler characteristic of $K$. The entries of
$h^{(c)}(K)$ can be expressed in terms of those of $h^{(sc)}(K)$ by
the equation
\begin{equation}\label{h^c synarthsei h^sc}
h^{(c)}_i(K)\; = \;\sum_{j=0}^{i-1}(-1)^{i+j-1}
h_j^{(sc)}(K)\,+\,(-1)^i 2^{d-1}, \quad 1\leq i\leq d.
\end{equation}
The cubical $h$-vectors of $K$ can also be expressed in terms of the
simplicial $h$-vectors of the links of the vertices of $K$ (these
links are simplicial complexes); see \cite[Theorem 9]{adin}.

\subsection{Cubical barycentric subdivision} \label{cbs-subsection} Let $K$ be a cubical
complex. The \emph{cubical barycentric subdivision}
$\mathrm{sd}_c(K)$ of $K$ is the polyhedral complex defined as
follows: Vertices of $\mathrm{sd}_c(K)$ are the barycenters of the
nonempty faces of $K$. Furthermore, to each closed interval $[F,G]$
in the poset $\mathcal{F}(K)\smallsetminus\{\varnothing\}$ of
nonempty faces of $K$ corresponds a face of $\mathrm{sd}_c(K)$. This
face is equal to the convex hull of the barycenters of all elements
of $[F,G]$; see Figure 1 for an example.

\vspace{0.8cm}

\includegraphics[bb= -30.000 0.000 100.000 120.000, scale=0.8]{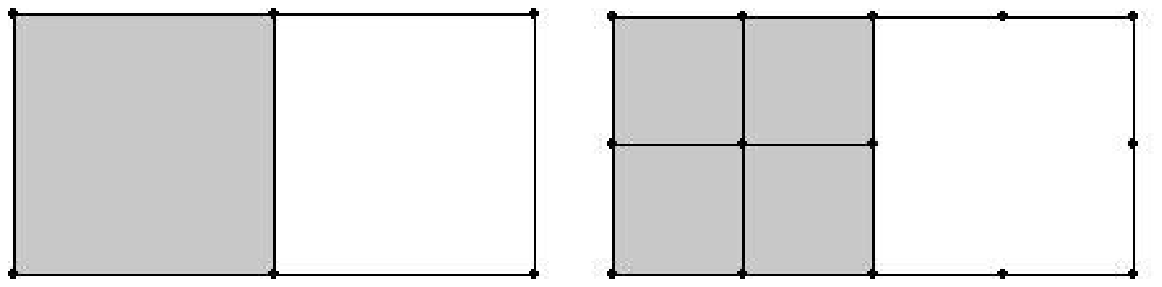}
\begin{center}
Figure 1: A cubical complex and its cubical barycentric subdivision
\end{center}
\vspace{0.7cm}

The face poset of $\mathrm{sd}_c(K)$ is isomorphic to the set of
closed intervals in the poset $\mathcal{F}(K)\smallsetminus
\{\varnothing\}$ of nonempty faces of $K$, partially ordered by
inclusion, with a minimum element adjointed. This fact implies that
$\mathrm{sd}_c(K)$ is also a cubical complex; its dimension is equal
to that of $K$.

\section{The short cubical $h$-vector} \label{section h^sc-vector}
This section studies the short cubical $h$-polynomial of
$\mathrm{sd}_c(K)$ with respect to its nonnegativity, symmetry and
real rootedness. The asymptotic behavior of the short cubical
$h$-polynomial of the $n^{\textrm{th}}$ iterated cubical barycentric
subdivision of $K$, as $n$ approaches infinity, is also explored.
Here and in the following section, $K$ denotes a $(d-1)$-dimensional
cubical complex and $\mathrm{sd}_c(K)$ stands for its cubical
barycentric subdivision.

To express the short cubical $h$-vector of $\mathrm{sd}_c(K)$ in
terms of that of $K$, we first need to determine the relationship
between the corresponding $f$-vectors.

\begin{prop} The $f$-vectors of $K$ and $\mathrm{sd}_c(K)$ are
related as follows:
\begin{equation}\label{relfvec}
f_i(\mathrm{sd}_c(K)) \; =\;  2^i \sum_{j=i}^{d-1} \binom{j}{i}
f_j(K), \quad i=0, \ldots, d-1.
\end{equation}
\end{prop}

\noindent \emph{Proof.} We recall that there is a one to one
correspondence between the set of $i$-dimensional faces of
$\mathrm{sd}_c(K)$ and the set of closed intervals $[x,y]$ of rank
$i$ in the face poset of $K$. To count these intervals, we note that
there are $f_j(K)$ ways to choose a face $y$ of $K$ of given
dimension $j$. Since every such face is a combinatorial $j$-cube,
there are $2^i\binom{j}{i}$ ways to choose a face $x$ of $y$ of
codimension $i$. As a result, there are $2^i\binom{j}{i}f_j(K)$
intervals $[x,y]$ of rank $i$ in the face poset of $K$ with
$\dim(y)=j$. Summing over all $j\geq i$, we obtain (\ref{relfvec}).
\hspace{\fill} $\square$

\begin{thm}\label{hsc}
The short cubical $h$-vectors of $K$ and $\mathrm{sd}_c(K)$ are
related as follows:
\begin{equation}\label{sxesh h^sc} h^{(sc)}_i(\mathrm{sd}_c(K)) \; = \;
\sum_{j=0}^{d-1} B(d,i,j)\, h^{(sc)}_j(K), \quad i=0, \ldots, d-1,
\end{equation}
where the coefficients $B(d,i,j)$ are nonnegative rational numbers,
determined by the generating function
\begin{equation}\label{genfunB}
\sum_{i=0}^{d-1}B(d,i,j)x^i \; =\; \frac{1}{2^{d-1}}\,
(3x+1)^j(x+3)^{d-1-j}\;.
\end{equation}
\end{thm}

\noindent \emph{Proof.} Equations (\ref{relfvec}) can be rewritten
as
\begin{equation}\label{relfpol}
f_{\mathrm{sd}_c(K)}(x)\;=\; f_K(1+2x),
\end{equation}
where $f_{K}(x)$ and $f_{\mathrm{sd}_c(K)}(x)$ are the
$f$-polynomials of $K$ and $\mathrm{sd}_c(K)$, respectively. Using
(\ref{syndesh polyn. short h me f}) and (\ref{syndesh polyn. f me
short h}), we can further rewrite (\ref{relfpol}) as
\begin{equation}\label{schpol4sd_c}
2^{d-1} h^{(sc)}_{\mathrm{sd}_c(K)}(x) \;=\; (x+3)^{d-1} h^{(sc)}_K
\left( \frac{3x+1}{x+3} \right),
\end{equation}
so that
\[
h^{(sc)}_{\mathrm{sd}_c(K)}(x) \; = \; \frac{1}{2^{d-1}}\;
\sum_{j=0}^{d-1} h_j^{(sc)}(K) (3x+1)^j (x+3)^{d-1-j}.
\]
The result follows by equating the coefficients of $x^i$ in the two
hand sides of this equation. \hspace{\fill} $\square$
\\

Recall that a vector $(a_0,a_1,\ldots,a_{n-1})\in\mathbb{R}^n$ is
said to be nonnegative (respectively, symmetric) if $a_i\geq 0$
(respectively, $a_i=a_{n-1-i}$) holds for $0\leq i\leq n-1$.

\begin{cor}\label{pos4short.cor}
If $K$ has nonnegative short cubical $h$-vector, then so does
$\mathrm{sd}_c(K)$.
\end{cor}

\noindent \emph{Proof.} Equation (\ref{genfunB}) implies that the
coefficients $B(d,i,j)$ are nonnegative for all $0\leq i,j \leq
d-1$. Thus the statement follows from (\ref{sxesh h^sc}).
\hspace{\fill} $\square$
\begin{cor}\label{symmetry for short.cor}
If $K$ has symmetric short cubical $h$-vector, then so does
$\mathrm{sd}_c(K)$.
\end{cor}

\noindent \emph{Proof.} Replacing $x$ by $1/x$ in (\ref{genfunB})
and multiplying by $x^{d-1}$, we find that
\[B(d,i,j)\; = \;B(d,d-1-i,d-1-j).\]
Assume that $h^{(sc)}(K)$ is symmetric, so that
$h^{(sc)}_{d-1-j}(K)=h^{(sc)}_j(K)$ holds for $0\leq j\leq d-1$.
Then
\begin{eqnarray*}
h^{(sc)}_{d-1-i}(\mathrm{sd}_c(K)) &=&
\sum_{j=0}^{d-1}B(d,d-1-i,j)\,h^{(sc)}_j(K)\\
&=&
\sum_{k=0}^{d-1}B(d,d-1-i,d-1-k)\,h^{(sc)}_{d-1-k}\\
&=& \sum_{k=0}^{d-1} B(d,i,k)\, h^{(sc)}_k \; =\;
h^{(sc)}_i(\mathrm{sd}_c(K))
\end{eqnarray*}
and hence $h^{(sc)}(\mathrm{sd}_c(K))$ is also symmetric.
\hspace{\fill} $\square$

\begin{cor}\label{realroot4hsc}
The short cubical $h$-polynomial of $K$ has only real roots if and
only if the same holds for $\mathrm{sd}_c(K)$.
\end{cor}

\noindent \emph{Proof.} This statement follows easily from
(\ref{schpol4sd_c}). The details are left to the reader.
\hspace{\fill} $\square$

\begin{ex}
{\rm For the boundary complex $K$ of the 3-dimensional cube we have
$h_K^{(sc)}(x) = 8 (1+x+x^2)$. This polynomial has positive
coefficients and two non-real complex roots. By Corollaries
\ref{pos4short.cor} and \ref{realroot4hsc}, so does the short
cubical $h$-polynomial of $\mathrm{sd}_c(K)$. This situation is in
contrast with what holds for barycentric subdivisions of simplicial
complexes (see \cite[Theorem 2]{brentiwelker}).}
\end{ex}

We denote by $\mathrm{sd}_c^n(K)$ the $n^{\textrm{th}}$ iterated
cubical barycentic subdivision of $K$, i.e. $\mathrm{sd}_c^0(K)=K$
and $\mathrm{sd}_c^n(K)=\mathrm{sd}_c(\mathrm{sd}_c^{n-1}(K))$ for
$n\geq1$. The short cubical $h$-polynomial of $\mathrm{sd}_c^n(K)$
has the following simple expression, in terms of the short cubical
$h$-polynomial of $K$.
\begin{prop}\label{propschpol4sd_cn}
The short cubical $h$-polynomials of $K$ and $\mathrm{sd}_c^n(K)$
are related as follows:
\begin{equation}\label{schpol4sd_cn} h^{(sc)}_{\mathrm{sd}_c^n(K)}(x) =
\left(\frac{(2^n-1)x+2^n+1}{2}\right)^{d-1}
h^{(sc)}_K\left(\frac{(2^n+1)x+2^n-1}{(2^n-1)x+2^n+1}\right).
\end{equation}
\end{prop}

\noindent \emph{Proof.} This follows from (\ref{schpol4sd_c}) by
induction on $n$. \hspace{\fill} $\square$
\\

The following statement implies that all complex roots of the short
cubical $h$-polynomial of $\mathrm{sd}_c^n(K)$ converge to $-1$ as
$n\rightarrow\infty$.

\begin{cor}\label{sch-polynomial 4 n-th subdivision}
We have
\[
\frac{1}{2^{n(d-1)}} \; h^{(sc)}_{\mathrm{sd}_c^n(K)}(x)
\;\rightarrow\; f_{d-1}(K)\; (x+1)^{d-1}
\]
coefficientwise, as $n\rightarrow\infty$. In particular, the short
cubical $h$-polynomial of $\mathrm{sd}_c^n(K)$ has positive and
unimodal coefficients for all large $n$.
\end{cor}
\noindent \emph{Proof.} The first statement follows from
(\ref{schpol4sd_cn}) and the fact that $h^{(sc)}_K(1)=2^{d-1}\,
f_{d-1}(K)$ by straightforward computations. The second statement
follows from the first and well-known unimodality properties of
binomial coefficients. \hspace{\fill} $\square$

\begin{rem}
{\rm We do not know of an example of cubical complex $K$ for which
$h^{(sc)}(K)$ is nonnegative and $h^{(sc)}(\mathrm{sd}_c(K))$ is not
unimodal.}
\end{rem}

\section{The cubical $h$-vector}\label{section h^c vector}

This section proves results on the cubical $h$-vector of
$\mathrm{sd}_c(K)$ analogous to those on the short cubical
$h$-vector in Section \ref{section h^sc-vector}.
\begin{thm}\label{thm h^c}
The cubical $h$-vectors of $K$ and $\mathrm{sd}_c(K)$ are related as
follows:
\begin{equation}\label{hc ths sd synarthsei hc}
h^{(c)}_i(\mathrm{sd}_c(K)) \; = \; \sum_{j=0}^d
C(d,i,j)h^{(c)}_j(K),
\end{equation}
where the coefficients $C(d,i,j)$ are nonnegative rational numbers,
determined by the generating function
\begin{eqnarray}\label{gen C(d,i,j) me x,y}
\sum_{j=0}^d \sum_{i=0}^d C(d,i,j) x^i y^j &=& \frac{1}{1+x} (1 +
x^{d+1}y^d) + \, \frac{x y}{2^{d-3}} \frac{(x+3)^{d-1} -
(3x+1)^{d-1} y^{d-1}}{x+3-(3x+1)y} \nonumber
\\ & &+\,\frac{x}{2^{d-1}(1+x)} ((x+3)^{d-1} +
(3x+1)^{d-1} y^d).
\end{eqnarray}
\end{thm}

\noindent \emph{Proof.} Since
$h^{(c)}_0(K)=h^{(c)}_0(\mathrm{sd}_c(K))=2^{d-1}$, equation
(\ref{hc ths sd synarthsei hc}) is valid for $i=0$ if we set
$C(d,0,0)=1$ and $C(d,0,j)=0$ for $1\leq j\leq d$. This agrees with
(\ref{gen C(d,i,j) me x,y}), since the right-hand side reduces to
the constant polynomial with value $1$ for $x=0$. Using (\ref{h^sc
synarthsei h^c}), (\ref{h^c synarthsei h^sc}) and (\ref{sxesh
h^sc}), we find that for $1\leq i \leq d$

\begin{eqnarray*}
h_i^{(c)}(\mathrm{sd}_c(K)) &=& \sum_{k=0}^{i-1}(-1)^{i+k-1}
\,h^{(sc)}_k(\mathrm{sd}_c(K)) \,+\, (-1)^i\,2^{d-1}\\
&=& \sum_{k=0}^{i-1} (-1)^{i+k-1} \sum_{j=0}^{d-1}B(d,k,j)\,h^{(sc)}_j(K)\, +\, (-1)^i\,2^{d-1}\\
&=& \sum_{k=0}^{i-1} (-1)^{i+k-1}
\sum_{j=0}^{d-1}B(d,k,j)(h^{(c)}_j(K) +
h^{(c)}_{j+1}(K))\, +\, (-1)^ih^{(c)}_0(K)\\
&=& \sum_{j=0}^d C(d,i,j) h^{(c)}_j(K),
\end{eqnarray*}
where
\begin{equation}
C(d,i,0) \; =\; \sum_{k=0}^{i-1} (-1)^{i+k-1} B(d, k, 0) + (-1)^i
\end{equation}
and
\begin{equation}\label{C(d,i,j)}
C(d,i,j)\; = \;\sum_{k=0}^{i-1}(-1)^{i+k-1}(B(d,k,j)+B(d,k,j-1))
\end{equation}
for $1\leq j\leq d$, under the convention that $B(d,k,d)=0$ for all
$k$. Using (\ref{genfunB}), we compute that for $1\leq j\leq d-1$

\begin{eqnarray*}
\sum_{i=0}^d C(d,i,j)x^i &=& \sum_{i=0}^d
\sum_{k=0}^{i-1} (-1)^{i+k-1} (B(d,k,j)+B(d,k,j-1))\,x^i\\
&&\\
&=& \sum_{k=0}^{d-1} (-1)^{k-1} (B(d,k,j)+B(d,k,j-1)) (-x)^{k+1} \,
\frac{1-(-x)^{d-k}}{1+x}\\ & & \\ &=& \frac{1}{1+x} \left(
x\sum_{k=0}^{d-1} B(d,k,j)x^{k} +
(-x)^{d+1} \sum_{k=0}^{d-1} B(d,k,j)(-1)^k \right.\\
& & + \left. x\, \sum_{k=0}^{d-1} B(d,k,j-1)x^{k} + (-x)^{d+1}
\sum_{k=0}^{d-1} B(d,k,j-1)(-1)^k \right)\\& & \\ &=&
\frac{1}{2^{d-1}(1+x)} (x(3x+1)^j (x+3)^{d-1-j} +
(-x)^{d+1}(-2)^j2^{d-1-j} \\ &&+\, x(3x+1)^{j-1} (x+3)^{d-j} +
(-x)^{d+1}(-2)^{j-1}2^{d-j})
\end{eqnarray*}
and hence
\begin{equation}\label{genCfor i,j=1,...,d-1}
\sum_{i=0}^d C(d,i,j)x^i\;=\;
\frac{x(3x+1)^{j-1}(x+3)^{d-1-j}}{2^{d-3}}.
\end{equation}
Similar computations yield
\begin{equation}\label{genCfor i,j=0}
\sum_{i=0}^d C(d,i,0)x^i \;=\; \frac{1}{1+x}
\left(\frac{x(x+3)^{d-1}}{2^{d-1}}+1\right)
\end{equation}
and
\begin{equation}\label{genCfor i,j=d} \sum_{i=0}^d C(d,i,d)x^i \; =\; \frac{1}{1+x}
\left(\frac{x(3x+1)^{d-1}}{2^{d-1}}+x^{d+1}\right).
\end{equation}
From these equations we can infer that the $C(d,i,j)$ are
nonnegative rational numbers for all $0\leq i,j\leq d$. Multiplying
the equations (\ref{genCfor i,j=1,...,d-1}), (\ref{genCfor i,j=0})
and (\ref{genCfor i,j=d}) by $y^j$ for $1\leq j\leq d-1$, $j=0$ and
$j=d$, respectively, and summing over all $j$ results in (\ref{gen
C(d,i,j) me x,y}).\hspace{\fill} $\square$

\begin{cor}\label{pos4long.cor}
If $K$ has nonnegative cubical $h$-vector, then so does
$\mathrm{sd}_c(K)$. \hspace{\fill} $\square$
\end{cor}

\begin{cor}\label{symmetry for long.cor}
If $K$ has symmetric cubical $h$-vector, then so does
$\mathrm{sd}_c(K)$.
\end{cor}

\noindent \emph{Proof.} As in the proof of Corollary \ref{symmetry
for short.cor}, it suffices to show that $C(d,d-i,d-j) = C(d,i,j)$
for all $0\leq i,j \leq d$. This follows from (\ref{gen C(d,i,j) me
x,y}) by replacing $x$ and $y$ by $1/x$ and $1/y$, respectively, and
multiplying by $x^{d}y^{d}$. \hspace{\fill} $\square$

\begin{rem}
{\rm We do not know of an example of cubical complex $K$ with
nonnegative cubical $h$-vector for which
$h^{(c)}_{\mathrm{sd}_c(K)}(x)$ is not real-rooted.}
\end{rem}

\begin{cor}\label{ch-polynomial 4 n-th subdivision}
For $d\geq 2$, we have
\begin{equation}
\frac{1}{2^{n(d-1)}}\; h^{(c)}_{\mathrm{sd}_c^n(K)}(x) \;
\rightarrow \; f_{d-1}(K)\, x\, (x+1)^{d-2}
\end{equation}
coefficientwise, as $n\rightarrow \infty$. In particular, if
$(-1)^{d-1}\tilde{\chi}(K)\geq 0$, then the cubical $h$-polynomial
of $\mathrm{sd}_c^n(K)$ has nonnegative and unimodal coefficients
for all large $n$.
\end{cor}

\noindent \emph{Proof.} Since $\mathrm{sd}_c(K)$ is a subdivision of
$K$, we have $\tilde{\chi}(\mathrm{sd}_c(K))=\tilde{\chi}(K)$ (this
equality also follows from (\ref{relfpol}) by setting $x=-1$). Thus,
by applying (\ref{rel short-long h-pol}) to $\mathrm{sd}_c^n(K)$, we
get
\begin{equation}\label{rel short-long h-pol 4sd^n}
(1+x) h^{(c)}_{\mathrm{sd}_c^n(K)}(x) \; = \; 2^{d-1} + x
h^{(sc)}_{\mathrm{sd}_c^n(K)}(x) + 2^{d-1} (-x)^{d+1}
\tilde{\chi}(K).
\end{equation}
Note that if $(-1)^{d-1}\tilde{\chi}(K)\geq 0$, then
$h_d^{(c)}(\mathrm{sd}_c^n(K))=(-2)^{d-1}\tilde{\chi}(K) \geq 0$.
The result follows from (\ref{rel short-long h-pol 4sd^n}) and
Corollary \ref{sch-polynomial 4 n-th subdivision}. \hspace{\fill}
$\square$

\end{document}